\documentclass[12pt,a4paper,oneside]{amsart}
\usepackage{amsfonts, amsmath, amssymb, amsthm, amscd}
\usepackage{anysize}
\usepackage{mathtools}
\marginsize{2cm}{2cm}{1.5cm}{1.5cm}

\usepackage[utf8]{inputenc}
\usepackage[english]{babel}
\usepackage{enumerate}

\usepackage{xcolor}
\usepackage{hyperref}
\hypersetup{
	colorlinks   = true, 
	urlcolor     = blue, 
	linkcolor    = blue, 
	citecolor    = red 
}
\newcommand{\Href}[2]{\hyperref[#2]{#1~\ref{#2}}}
\usepackage{anysize}
\marginsize{2cm}{2cm}{1.5cm}{1.5cm}

\makeatletter
\@namedef{subjclassname@2020}{\textup{2020} Mathematics Subject Classification}
\makeatother

\newtheorem{thm}{Theorem}

\newtheorem{cor}{Corollary}[section]

\newtheorem{claim}{Claim}[section]

\theoremstyle{definition}

\newcommand{\st}{:\;}

\def\R{{\mathbb R}}%
\renewcommand{\Re}{\mathbb{R}}
\newcommand{\ball}[2][d]{\mathbf{B}^{#1}\left(#2\right)}

\providecommand{\parenth}[1]{\left(#1\right)}%
\newcommand{\iprod}[2]{\left\langle#1,#2\right\rangle}%
\newcommand{\conv}[1]{\mathrm{conv}\left(#1\right)}
\newcommand{\sconv}[1]{\mathrm{Sconv}\left(#1\right)}
\newcommand{\Sed}{\mathbb{S}^{d}}
\newcommand{\Sedn}{\mathbb{S}^{d}_N}

\newcommand{\scap}[1]{\mathrm{Cap}\!\left(#1\right)}
\newcommand{\Red}{\Re^{d}}
\newcommand{\Redp}{\Re^{d+1}}
\newcommand{\polar}[1]{#1^{\circ}}
\def\spolarsign{+}%
\newcommand{\spolar}[1]{{#1}^{\spolarsign}}%

\newcommand{\setdef}[1]{\left\{#1\right\}}
\newcommand{\arccot}[1]{\mathrm{arccot}\!\left(#1\right)}

\title{Quantitative Steinitz theorem: A spherical version}
\author{Grigory Ivanov\address{Grigory Ivanov: 
Institute of Science and Technology Austria (IST Austria), 
Klosterneuburg, 3400, Austria}
\email{grimivanov@gmail.com}
\and
M\'arton Nasz\'odi\address{M\'arton Nasz\'odi:
Alfr\'ed R\'enyi Inst. of Mathematics HUN-REN and
Dept. of Geometry, Lor\'and E\"otv\"os University, Budapest}
\email{marton.naszodi@renyi.hu}
}

\thanks{
M.N. was supported by the J\'anos Bolyai Scholarship of the Hungarian Academy of Sciences as well as the National Research, Development and Innovation Fund (NRDI) grants K119670, K131529 and K147544, and the \'UNKP-22-5 New National Excellence Program of the Ministry for Innovation and Technology from the source of the NRDI.
}


\begin{document}
\begin{abstract}
Steinitz's theorem states that if the origin belongs to the interior of the convex hull of a set $Q \subset \mathbb{R}^d$, then there are at most $2d$ points $Q^\prime$ of $Q$ whose convex hull contains the origin in the interior. 
B\'ar\'any, Katchalski and Pach gave a quantitative version whereby the radius of the ball contained in the convex hull of $Q^\prime$ is bounded from below.
In the present note, we show that a Euclidean result of this kind implies a corresponding spherical version.
\end{abstract}
\maketitle

\section{Introduction}

A fundamental fact in convexity discovered by Steinitz states that the interior of the convex hull of a subset $Q$ of $\Red$ equals the union of the interiors of the convex hulls of at most $2d$ points of $Q$. 
A quantitative version was shown by B\'ar\'any, Katchalski and Pach \cite{barany1982quantitative, barany1984helly}. We will call it the Euclidean Quantitative Steinitz Theorem or, QST in short.

\begin{thm}[Euclidean QST]\label{thm:QST_monochromatic}
Let $Q$ be a subset of $\Red$  whose convex hull contains the Euclidean unit ball $\ball{o,1}$ centered at the origin $o$.
Then there exists set $Q^{\prime}$ of at most $2d$ points of $Q$ that satisfies
\[
 \ball{o,r} \subset \conv{Q^\prime}
\]  
with some $r>0$ depending only on $d$.
\end{thm}

We will use $r(d)$ to denote the largest $r$ that makes the conclusion of \Href{Theorem}{thm:QST_monochromatic} true. Clearly, $r(d) \leq 1.$

In \cite{barany1982quantitative}, the lower bound $r(d)>d^{-2d}$ is presented, $r(d)>cd^{-1/2}$ is conjectured with a universal constant $c>0$. The first polynomial lower bound, $r(d) > \frac{1}{6d^2}$ was proved in \cite{ivanov2022steinitz} by the authors. The upper bound $r(d)\leq\frac{1}{2\sqrt{d}}$ was shown in \cite[Theorem 2]{ivanov2022steinitz} and thus, $r(d)$ tends to zero as $d$ tends to infinity.

In the present note, we establish a spherical (or, cone) version of \Href{Theorem}{thm:QST_monochromatic}. We use $\Sed=\setdef{u\in\Redp\st \iprod{u}{u}=1}$ to denote the unit sphere in $\Redp$, and $e_{d+1}$ to denote the last element of the standard basis.
We define the \emph{spherical cap} with center $v\in\Sed$ and spherical radius $\rho\in[0,\pi]$ as
\[
\scap{v,\rho} = \setdef{u \in \Sed \st \iprod{u}{v} \geq \cos{\rho}}.
\]
A set $K\subset\Sed$ is called \emph{spherically convex}, if it is either $\Sed$, or is contained in an open hemisphere, and for any two points $u$ and $v$ of $K$, the shorter great circular arc connecting $u$ and $v$ is contained in $K$. The \emph{spherical convex hull} of a subset $A$ of $\Sed$, denoted by $\sconv{A}$, is defined accordingly. It is known to exist and to be unique for any $A\subseteq\Sed$.

\begin{thm}[Spherical QST]\label{thm:QST_monochromaticS}
Let $C$ be a subset of $\Sed$, with $d\geq2$, whose spherical convex hull contains the cap $\scap{e_{d+1},\rho}$ for some $\rho\in(0,\pi/2)$.
Then there exists set $C^\prime$ of at most $2d$ points of $C$ that satisfies
\[
\scap{e_{d+1},\gamma\rho} \subset \sconv{C^\prime}
\]  
with some $\gamma>0$ depending only on $d$.
\end{thm}

We will use $\gamma(d)$ to denote the largest $\gamma$ that makes the conclusion of \Href{Theorem}{thm:QST_monochromaticS} true.
Our main result is as follows.

\begin{thm}[Spherical QST follows from Euclidean QST]\label{thm:SphericalQST_from_Euclidean_QST} With the notation above,
\[
r(d) \geq \gamma(d)\geq \frac{r(d)}{2} 
\]
for any dimension $d\geq2$.
\end{thm}

Clearly, \Href{Theorem}{thm:SphericalQST_from_Euclidean_QST} implies \Href{Theorem}{thm:QST_monochromaticS}. Moreover, combined with the lower bound on $r(d)$ in the paragraph following \Href{Theorem}{thm:QST_monochromatic}, we obtain the following explicit bound.

\begin{cor}[Spherical QST]\label{cor:Spherical_QST_monochromatic}
Let $C$ be a subset of $\Sed$, with $d\geq2$, whose spherical convex hull contains the cap $\scap{e_{d+1},\rho}$ for some $\rho \in (0, \pi/2)$. 
Then there exists set $C^\prime$ of at most $2d$ points of $C$ that satisfies
\[
\scap{e_{d+1},\rho\frac{r(d)}{2}} \subset \sconv{C^\prime}.
\]
\end{cor}

While a non-quantitative spherical version of Steinitz's non-quantitative Euclidean theorem easily follows, the proof of our quantitative spherical result requires some additional ideas described in \Href{Section}{sec:SphericalQST_from_Euclidean_QS}

\section{Proof of \texorpdfstring{\Href{Theorem}{thm:SphericalQST_from_Euclidean_QST}}{{Theorem}~\ref{thm:SphericalQST_from_Euclidean_QST}}}\label{sec:SphericalQST_from_Euclidean_QS}
We will use the following  two simple and purely technical observations.
We postpone their proofs to the next section.
\begin{claim}\label{claim:tan_pi_4}
 For any $t>0$, we have $t \leq \tan t$.
For any $t \in \parenth{0, \frac{\pi}{4}}$, we have $\tan t  \leq  2t$.
\end{claim}
\begin{claim}\label{claim:tanestimate}
Assume that $ r \in (0,1]$  and $\rho \in \parenth{0, \frac{\pi}{2}}.$
Then 
\[
\frac{\pi}{2} - \arctan \parenth{\frac{\cot \rho}{r}} \geq 
\frac{r\rho}{2}.
\]
\end{claim}

A straight-forward combination of Euclidean polarity with \Href{Theorem}{thm:QST_monochromatic} yields the following Quantitative Helly Theorem (QHT), which we will use. We omit the proof.
We recall that the \emph{polar} of the set $S \subset \Red$ is defined by
\[
\polar{S} = \left\{x \in \Red \st \iprod{x}{s} \leq 1 \quad \text{for all} \quad s \in S
\right\}.
\]

\begin{thm}[Euclidean QHT]\label{thm:QHT}
Let $L$ be a subset of $\Red$  with $\polar{L}$ contained in $\ball{o,1}$.
Then there exists a set $L^{\prime}$ of at most $2d$ points of $L$ that satisfies
\[
 \ball{o,\frac{1}{r(d)}} \supset \polar{(L^\prime)},
\]
where $r(d)$ is the quantity defined after \Href{Theorem}{thm:QST_monochromatic}.
\end{thm}

We denote the Northern hyperplane by $H_N = \setdef{x \in \Redp \st \iprod{x}{e_{d+1}} = 1}$ and the open Northern hemisphere by $\Sedn = \{x \in \Redp \st \iprod{x}{e_{d+1} } >  0 \}$.
We identify $H_N$ with $\Red$ inheriting the metric of $\Redp$ and setting the point $e_{d+1}$ as the origin. 
We use $P_N$ to denote the central projection from the origin of $\Redp$ onto $H_N$, that is, for any $(x_1,\dots,x_{d+1})\in\Redp$ with $x_{d+1}\neq0$, we set $P_N\big((x_1,\dots,x_{d+1})\big)=\left(\frac{x_1}{x_{d+1}},\frac{x_2}{x_{d+1}}, \dots, \frac{x_d}{x_{d+1}}\right)$. 

First, we discuss the easy case, when $C$ is contained in $\Sedn$.
In this case the result directly follows from its Euclidean counterpart.
 We set $Q = P_N (C)$ and observe that $\conv{Q}\supseteq P_N(\scap{e_{d+1},\rho})=\ball{o,\tan \rho}$.
By the definition of $r(d)$, there is a subset $C^\prime$ of $C$
of size at most $2d$ such that $\conv{P_N (C^\prime)}$ contains the ball
$\ball{o,r(d) \tan{\rho}}$. Equivalently, 
$\sconv{C^\prime}$ contains the cap $\scap{e_{d+1}, \arctan \parenth{r(d) \tan{\rho}}}$.
Since $r(d) \leq 1,$ one gets $\frac{r(d)}{2} \cdot \rho < \frac{\pi}{2}.$  
Thus, by \Href{Claim}{claim:tan_pi_4}, the inequality
\[
\tan \parenth{\frac{r(d)}{2} \cdot \rho} \leq  r(d) \cdot \rho \leq r(d) \tan{\rho}
\]
holds for all $\rho \in (0, \pi/2)$ and all positive integer $d$. 
That is,  for all $\rho \in (0, \pi/2)$ and all positive integer $d,$ 
$\sconv{C^\prime}$ contains the cap $\scap{e_{d+1}, \frac{r(d)}{2} \cdot \rho}$.
However, since $\lim\limits_{t \to 0} \frac{\arctan{t}}{t} = \lim\limits_{t \to 0} \frac{\tan{t}}{t}=1,$ this case ensures that $\gamma(d) \leq r(d).$ 

Next, we discuss the general case, that is, when $C$ is not necessarily contained in the open Northern hemisphere. 

 We will use the following notation. For a point $c\in\Redp$ and a set $C\subset\Redp$, we set
\[
 \spolar{c}=\setdef{x\in\Redp\st \iprod{x}{c}>0}, \text{ and }
 \spolar{C}=\bigcap_{c\in C}\spolar{c}.
\]

Observe that for a set $C\subseteq\Sed$ we have that $\spolar{C}$ is empty if and only if, $C$ is not contained in \emph{any} open hemisphere. This, in turn, is equivalent to having $o\in\conv{C}$ in $\Redp$. If $C$ is such, then by Carath\'eodory's theorem, there is a $C^{\prime}\subseteq C$ of size at most $d+2$ with $o\in\conv{C^{\prime}}$. Which yields that $\sconv{C^{\prime}}=\Sed$, and there is nothing to prove.

Thus, we will assume that $K:=\spolar{C}\cap\Sed$ is not empty.
Then $K\subset \scap{e_{d+1},\frac{\pi}{2}-\rho}$, and thus, by projection, we have in $H_N$
\[P_N(K)\subset P_N\left(\scap{e_{d+1},\frac{\pi}{2}-\rho}\right)=\ball{e_{d+1},\cot(\rho)}.\]

Since $K\subset \Sedn$, we can write
\[ 
P_N(K)=P_N\left(\bigcap_{c\in C}\spolar{c}\right)=
\bigcap_{c\in C}P_N\left(\spolar{c}\cap\Sedn\right),
\]
where on the right, we see an intersection of half-spaces in $\Red$. Applying \Href{Theorem}{thm:QHT}, we obtain a set $C^{\prime}\subset C$ of size at most $2d$ with
\[ 
\ball{e_{d+1},\frac{\cot(\rho)}{r(d)}}\supset
\bigcap_{c\in C^{\prime}}P_N\left(\spolar{c}\cap\Sedn\right),
\]
which yields
\[ 
\scap{e_{d+1},\arctan \left(\frac{\cot(\rho)}{r(d)}\right)}\supset
\bigcap_{c\in C^{\prime}}\left(\spolar{c}\cap\Sed\right).
\]
By polarity on $\Sed$, we get
\[ 
\scap{e_{d+1},\frac{\pi}{2} - \arctan\left(\frac{\cot(\rho)}{r(d)}\right)}\subset
\sconv{C^{\prime}}.
\]

By \Href{Claim}{claim:tanestimate},
\[
 \frac{\pi}{2}-\arctan\left(\frac{\cot(\rho)}{r(d)}\right)\geq\frac{r(d)\rho}{2},
\]
and \Href{Theorem}{thm:SphericalQST_from_Euclidean_QST} follows.

\section{Proofs of the technical tools}
\begin{proof}[Proof of \Href{Claim}{claim:tan_pi_4}]
The first inequality is a standard exercise in calculus.
To prove the second one, consider $f \colon \left[0, \frac{\pi}{2}\right) \to \R$ given by
$f(t) = \frac{\tan t}{t}$ for $t \in \parenth{0, \frac{\pi}{2}}$  and $f(0)=1.$
By direct computations,
\[
f^\prime(t) = \frac{(\tan t)^\prime \cdot t - \tan t}{t^2} = 
\frac{\frac{t}{\cos^2 t} - \tan t}{t^2} = 
\frac{t - \sin t \cos t}{t^2 \cos^2 t} = \frac{t - \sin (2t)/2}{t^2 \cos^2 t}.
\]
Since $t > \sin t$ for a positive $t,$ we conclude that $f^\prime$ is positive on 
$(0, \pi/2).$ Thus, $f$ is strictly increasing on $\left[0, \frac{\pi}{2}\right).$
Consequently,  for any $t \in \parenth{0, \frac{\pi}{4}},$ 
\[
\tan t  \leq \frac{\tan (\pi/4)}{\pi/4} \cdot t \leq 2t.
\]
\end{proof}



\begin{proof}[Proof of \Href{Claim}{claim:tanestimate}]
Since $r > 0$ and $\rho \in  \parenth{0, \frac{\pi}{2}}$, we have the following sequence of equivalent inequalities:
\[
\frac{\pi}{2} - \arctan{\frac{\cot \rho}{r}} \geq \frac{r}{2} \cdot \rho   
\quad \Leftrightarrow \quad
\arccot{\frac{\cot \rho}{r}}
\geq \frac{r}{2} \cdot \rho 
\quad \Leftrightarrow \quad \frac{\cot \rho}{r}  \leq \cot\parenth{\frac{r}{2} \cdot \rho}
\quad \Leftrightarrow \quad 
\]
\[
\tan \parenth{\frac{r}{2} \cdot \rho} \leq r \tan {\rho}.
\]
Since $r \in (0,1],$ one gets $\frac{r}{2} \cdot \rho \in \parenth{0, \frac{\pi}{4}}.$
By Claim \ref{claim:tan_pi_4},
\[
\tan \parenth{\frac{r}{2} \cdot \rho} \leq r \rho \leq r \tan \rho.
\]
The claim follows.
\end{proof}

\bibliographystyle{amsalpha}
\bibliography{uvolit}
\end{document}